\documentclass{amsart}

\usepackage{amsthm, amssymb, amsfonts, latexsym}
\usepackage{hyperref}

\hypersetup{
pdftitle = {The vanishing ideal of a finite set of closed points in affine space},
pdfsubject = {Fakultaet fuer Mathematik, Universitaet Bielefeld, Germany},
pdfauthor = {Mathias Lederer},
pdfkeywords = {Gr\"obner bases, vanishing ideal, zero-dimensional ideal, 
Buchberger--M\"oller algorithm, lexicographic ordering, Lederer},
}

\newcommand{\N}{\,\mathbb{N}}
\newcommand{\A}{\,\mathbb{A}}
\newcommand{\D}{\mathcal{D}}
\newcommand{\Oh}{\mathcal{O}}

{\theoremstyle{definition} 
\newtheorem*{dfn}{Definition}}
\newtheorem{pro}{Proposition}
\newtheorem{lmm}{Lemma}
\newtheorem{thm}{Theorem}
\newtheorem*{cor}{Corollary}
\newtheorem*{ass}{Assumption}

\begin{document}

\title{The vanishing ideal of a finite set of closed points in affine space}
\author{Mathias Lederer}
\address{Fakult\"at f\"{u}r Mathematik, Universit\"{a}t Bielefeld, Bielefeld, Germany}
\email{mlederer@math.uni-bielefeld.de}
\subjclass{13P10, 14Q99, 14Q20, 14R10}
\date{April, 2006}
\keywords{Gr\"obner bases, vanishing ideal, zero-dimensional ideal, Buchberger--M\"oller algorithm, lexicographic ordering}

\maketitle

\begin{abstract} Given a finite set of closed rational points of affine space over a field, 
  we give a Gr\"obner basis for the lexicographic ordering of the ideal of polynomials 
  which vanish at all given points.   
  Our method is an alternative to the Buchberger--M\"oller algorithm, but in contrast to that, 
  we determine the set of leading terms of the ideal without solving any linear equation
  but by induction over the dimension of affine space. 
  The elements of the  Gr\"obner basis are also computed by induction over the dimension, 
  using one-dimensional interpolation of coefficients of certain polynomials.
\end{abstract}

\section{Introduction}\label{intro}

Let $k$ be a field. Consider the affine space $\A^n$ over $k$. Suppose we are given a finite set $A$ of closed
$k$-rational points of $\A^n$, i.e. each $a\in\A^n$ is given by coordinates $a=(a_{1},\ldots,a_{n})\in k^n$. 
Our aim is to find a Gr\"obner basis of the ideal 
\begin{equation*}
  I(A)=\{f\in k[X];\forall a\in A:\,f(a)=0\}
\end{equation*}
in $k[X]$, where we write $X=(X_{1},\ldots,X_{n})$ for brevity's sake. We will use the lexicographical ordering on $k[X]$, 
where $X_{1}<X_{2}<\ldots<X_{n}$ and give our  Gr\"obner basis solely for this particular ordering. 

There exists an algorithm that provides a complete solution to this problem -- to wit, the Buchberger--M\"oller algorithm. 
It first appeared in \cite{4} and was subsequently generalised in \cite{1} to apply to 
$k[X]$-modules and $k[X]$-submodules instead of $k[X]$ and ideals within. 
The Buch\-berger--M\"oller algorithm treats the problem of finding a Gr\"obner basis of $I(A)$ in a more general way than the
present article does, since already the original article \cite{4} makes no restriction on the term ordering on $k[X]$ 
for which the Gr\"obner basis is constructed. However, our answer for the special case of lexicographical ordering will be 
in a way more transparent than what the Buchberger--M\"oller algorithm does. In particular, 
we will explicitly know the set of leading terms of elements of $I$ from the relative position of the elements of $A$. 
Since in its complete form our construction is rather involved, it is recommendable to first illustrate the idea of the method 
by looking at a few examples. 

First, let us take $A=\{(1,0), (1,2), (3,1), (3,4)\}$. It is easy to write down one element of $I(A)$:
\begin{equation*}
  f_{1}=(X_{1}-1)(X_{1}-3)=X_{1}^2-4X_{1}+3\in I(A)\,.
\end{equation*}
In fact, for writing down this polynomial, we project the elements of $A$ to $\A^1$ by means of 
$p_{1}:(a_{1},a_{2})\mapsto a_{1}$ and then consider the ideal $I(p_{1}(A))$ in $k[X_{1}]$, 
whose Gr\"obner basis is trivial to compute.

Next, we might also try this for the projection $p_{2}$ instead of $p_{1}$. But this would substantially change the situation,
since $\#p_{2}(A)=4$, whereas $\#p_{1}(A)=2$. Here is a better idea. The two polynomials 
\begin{equation*}
  \begin{split}
    g&=X_{2}(X_{2}-2)=X_{2}^2-2X_{2}\text{ and }\\
    h&=(X_{2}-1)(X_{2}-4)=X_{2}^2-5X_{2}+4
  \end{split}
\end{equation*}
do not lie in $I(A)$, but at least $g$ vanishes on $p_{1}^{-1}(1)\cap A$ and $h$ on $p_{1}^{-1}(3)\cap A$. 
Therefore, let us modify the coefficients of $g$ and $h$ in such a way that the result, call it $f_{2}$, will vanish at all elements of $A$. 
This can also be done by applying a familiar technique from the one-dimensional case. 
We simply replace each coefficient of $g\in k[X_{2}]$, respectively of $h\in k[X_{2}]$, by the polynomial in $k[X_{1}]$
that interpolates the coefficient of $g$ and the corresponding coefficient of $h$. 
In other words, we use the characteristic polynomials
\begin{equation*}
  \begin{split}
    \chi_{1}&=\frac{X_{1}-3}{1-3}=-\frac{1}{2}(X_{1}-3)\text{ and }\\
    \chi_{3}&=\frac{X_{1}-1}{3-1}=\frac{1}{2}(X_{1}-1)
  \end{split}
\end{equation*}
of $1\in\{1,3\}$ and $3\in\{1,3\}$, respectively, to define
\begin{equation*}
  f_{2}=\chi_{1}g+\chi_{2}h\in k[X_{1},X_{2}]\,.
\end{equation*}
Then $f_{2}$ clearly lies in $I(A)$. Since $\chi_{1}+\chi_{2}=1$, the leading term of $f_{2}$ is $X_{2}^2$. 
The lower terms of $f_{2}$ are $k$-multiples of $X_{1}X_{2}$, $X_{2}$, $X_{1}$ and $1$, respectively. 
The leading term of $f_{i}$ in particular divides none of the nonleading terms of $f_{i}$, for $i,j\in\{1,2\}$.
Therefore, $(f_{1},f_{2})$ is a Gr\"obner basis of $I(A)$. (This reasoning is standard in the theory of Gr\"obner bases 
(\cite{2}, \cite{3}) and will henceforth be used without explicit mention.)

\begin{center}
\begin{figure}[ht]
  \begin{picture}(270,110)
    \put(30,10){\line(1,0){90}}
    \put(30,10){\line(0,1){90}}
    \multiput(30,10)(20,0){5}{\line(0,-1){7}}
    \multiput(30,10)(0,20){5}{\line(-1,0){7}}
    \put(50,10){\circle*{3}}
    \put(50,50){\circle*{3}}
    \put(90,30){\circle*{3}}
    \put(90,90){\circle*{3}}
    \put(150,10){\line(1,0){90}}
    \put(150,10){\line(0,1){90}}
    \multiput(150,10)(20,0){5}{\line(0,-1){7}}
    \multiput(150,10)(0,20){5}{\line(-1,0){7}}
    \put(150,10){\circle{3}}
    \put(150,30){\circle{3}}
    \put(170,10){\circle{3}}
    \put(170,30){\circle{3}}
    \multiput(150,50)(0,20){3}{\multiput(0,0)(20,0){2}{\circle*{3}}}
    \multiput(190,10)(0,20){5}{\multiput(0,0)(20,0){3}{\circle*{3}}}
  \end{picture}
\caption{The elements of $A$ and the exponents of $I(A)$}
\label{ex1}
\end{figure}
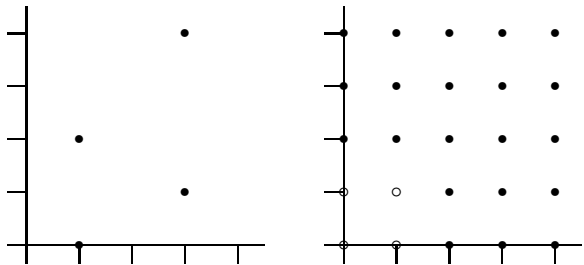
\end{center}

The left picture in Figure \ref{ex1} shows $A\subseteq\A^2$. The right picture shows those elements of
$\N_{0}^2$ that occur as exponents of leading terms of $I(A)$. Each of these elements is marked by a solid circle.
Note that (as an additive submonoid of $\N_{0}^2$) this set is spanned by $(2,0)$ and $(0,2)$, the exponents of 
$f_{1}$ and $f_{2}$ respectively.

At this point a comment on the set in the right picture of Figure \ref{ex1} may be in order. 
What we will be working with is actually not the set of those elements of $\N_{0}^2$ that occur as 
exponents of leading terms of elements of $I(A)$ but rather its complement in $\N_{0}^n$, call it $D(A)$.
In the above picture, the elements of $D(A)$ are marked by blank circles. 
We have built up the polynomial $f_{2}$ by looking at $p_{1}^{-1}(1)\cap A$ and $p_{1}^{-1}(3)\cap A$. 
We understand these two sets as subsets of $\A^1$ by means of the projection $p_{2}: (a_{1},a_{2})\mapsto a_{2}$. 
This leads to the subsets $D(p_{1}^{-1}(1)\cap A)$ and $D(p_{1}^{-1}(3)\cap A)$ of $\N_{0}$
(analogously defined as $D(A)$). In our example, it becomes evident from $g$ and $h$, respectively, 
that $D(p_{1}^{-1}(1)\cap A)=D(p_{1}^{-1}(3)\cap A)=\{0,1\}$.
One key result of the present article is that $D(A)$ is built up from the two blocks $D(p_{1}^{-1}(1)\cap A)$ 
and $D(p_{1}^{-1}(3)\cap A)$ in a quite intuitive way, as shown in Figure \ref{add1}. 
This will be given precise definition in Section \ref{gamma}. 

\begin{center}
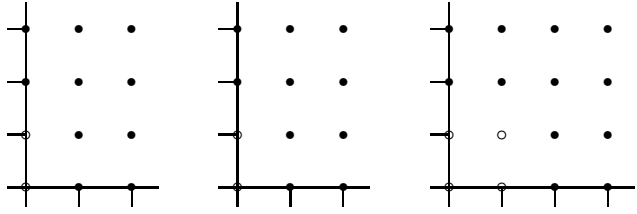
\begin{figure}[ht]
  \begin{picture}(290,90)
    \put(30,10){\line(1,0){50}}
    \put(30,10){\line(0,1){70}}
    \multiput(30,10)(20,0){3}{\line(0,-1){7}}
    \multiput(30,10)(0,20){4}{\line(-1,0){7}}
    \multiput(30,10)(0,20){2}{\circle{3}}
    \multiput(30,50)(0,20){2}{\circle*{3}}
    \multiput(50,10)(20,0){2}{\multiput(0,0)(0,20){4}{\circle*{3}}}
    \put(110,10){\line(1,0){50}}
    \put(110,10){\line(0,1){70}}
    \multiput(110,10)(20,0){3}{\line(0,-1){7}}
    \multiput(110,10)(0,20){4}{\line(-1,0){7}}
    \multiput(110,10)(20,0){3}{\line(0,-1){7}}
    \multiput(110,10)(0,20){4}{\line(-1,0){7}}
    \multiput(110,10)(0,20){2}{\circle{3}}
    \multiput(110,50)(0,20){2}{\circle*{3}}
    \multiput(130,10)(20,0){2}{\multiput(0,0)(0,20){4}{\circle*{3}}}
    \put(190,10){\line(1,0){70}}
    \put(190,10){\line(0,1){70}}
    \multiput(190,10)(20,0){4}{\line(0,-1){7}}
    \multiput(190,10)(0,20){4}{\line(-1,0){7}}
    \multiput(190,10)(20,0){2}{\multiput(0,0)(0,20){2}{\circle{3}}}
    \multiput(190,50)(20,0){2}{\multiput(0,0)(0,20){2}{\circle*{3}}}
    \multiput(230,10)(20,0){2}{\multiput(0,0)(0,20){4}{\circle*{3}}}
  \end{picture}
\caption{$D(p_{1}^{-1}(1)\cap A)$ and $D(p_{1}^{-1}(3)\cap A)$ together form $D(A)$}
\label{add1}
\end{figure}
\end{center}

Let us consider a second example. Take $A^\prime=\{(1,0), (1,2), (2,3), (3,1), (3,4)\}=A\cup\{(2,3)\}$. 
Again, the first element of $I(A^\prime)$ is easy to write.
\begin{equation*}
  f_{1}=(X_{1}-1)(X_{1}-2)(X_{1}-3)=X_{1}^3-6X_{1}^2+11X_{1}-6\in I(A^\prime)\,.
\end{equation*}
For imitating the construction of $f_{2}$, we first take the three polynomials 
\begin{equation*}
  \begin{split}
    g&=X_{2}(X_{2}-2)=X_{2}^2-2X_{2}\,,\\
    h&=(X_{2}-1)(X_{2}-4)=X_{2}^2-5X_{2}+4\text{ and }\\
    i&=(X_{2}-3)(X_{2}+3)=X_{2}^2-9\,,
  \end{split}
\end{equation*}
where $g$ and $h$ vanish on the same subsets of $A^\prime$ as before and $i$ vanishes on $p_{1}^{-1}(2)\cap A^\prime$. 
(We will presently see why $i=X_{2}-3$ would not be a good choice.) Now we need
\begin{equation*}
  \begin{split}
    \chi_{1}&=\frac{(X_{1}-2)(X_{1}-3)}{(1-2)(1-3)}=\frac{1}{2}(X_{1}^2-5X_{1}+6)\,,\\
    \chi_{2}&=\frac{(X_{1}-1)(X_{1}-3)}{(2-1)(2-3)}=-(X_{1}^2-4X_{1}+3)\text{ and }\\
    \chi_{3}&=\frac{(X_{1}-1)(X_{1}-2)}{(3-1)(3-2)}=\frac{1}{2}(X_{1}^2-3X_{1}+2)\,,
  \end{split}
\end{equation*}
the characteristic polynomials of $1,2$ and $3\in\{1,2,3\}$, respectively. We define
\begin{equation*}
  f_{2}=\chi_{1}g+\chi_{3}h+\chi_{2}i\,.
\end{equation*}
The leading term of $f_{2}$ is $X_{2}^2$, and the lower terms of $f_{2}$ are $k$-multiples of 
$X_{1}^2X_{2}$, $X_{1}X_{2}$, $X_{2}$, $X_{1}^2$, $X_{1}$ and $1$, respectively. 
As before, the pair $(f_{1},f_{2})$ is a Gr\"obner basis. But the dimension of $k[X]/(f_{1},f_{2})$ 
as a $k$-vector space is $6$, whereas the dimension of $k[X]/I(A)$ is $5$ 
(by the Chinese Remainder Theorem, since $\#A^\prime=5$).

Therefore, the ideal $(f_{1},f_{2})$ is bigger than $I(A^\prime)$. The reason for this should appear from the polynomial $i$
itself. In fact, it is necessary to have a polynomial $i$ that vanishes 
on $p_{1}^{-1}(2)\cap A^\prime$, whose leading term equals $X_{2}^2$. 
Only in this way can we guarantee that the leading term of $f_{2}=\chi_{1}f+\chi_{3}g+\chi_{2}i$ is $X_{2}^2$. 
We defined $i=(X_{2}-3)(X_{2}+3)$, which is of the desired shape -- but unfortunately, it vanishes not only at $3$. 
Its other zero is $-3$, so $(f_{1},f_{2})$ is a Gr\"obner basis of $I(A^{\prime\prime})$, 
where $A^{\prime\prime}=A\cup\{(2,-3)\}$. Defining $i=(X_{2}-3)^2$ would not make things better, 
since the dimension of $k[X]/(f_{1},f_{2})$ is also $6$ when this input is used.
(The reader will understand why we have taken $i=(X_{2}-3)(X_{2}+3)$ and nothing else of the same kind 
after Section \ref{almost}.)

The way out goes as follows: Set 
\begin{equation*}
  f_{3}=(X_{1}-1)(X_{1}-3)(X_{2}-3)\,.
\end{equation*}
This polynomial also lies in $I(A^\prime)$. The leading term of $f_{3}$ is $X_{1}^2X_{2}$ 
and its lower terms are $k$-multiples of $X_{1}X_{2}$, $X_{2}$, $X_{1}^2$, $X_{1}$ and $1$, 
respectively. Therefore, the linear combination $f_{2}-cf_{3}$, where $c=4$ is the coefficient of $X_{1}^2X_{2}$ 
in $f_{2}$, lies in $I(A^\prime)$ as well, but this polynomial has $X_{2}^2$ as leading term and 
$k$-multiples of $X_{1}X_{2}$, $X_{2}$, $X_{1}^2$, $X_{1}$ and $1$, respectively, as lower terms. 
Therefore, the $k$-dimension of $k[X]/(f_{1},f_{2},f_{3})$ is 5, 
hence $(f_{1},f_{2}-cf_{3},f_{3})$ is a Gr\"obner basis of $I(A^\prime)$. 

\begin{center}
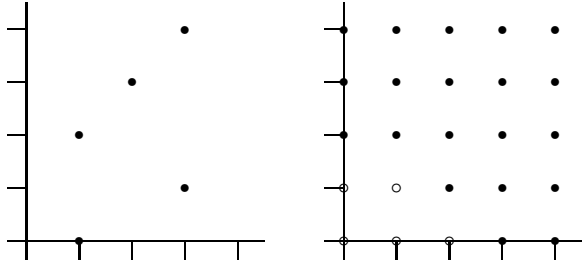
\begin{figure}[ht]
  \begin{picture}(270,110)
    \put(30,10){\line(1,0){90}}
    \put(30,10){\line(0,1){90}}
    \multiput(30,10)(20,0){5}{\line(0,-1){7}}
    \multiput(30,10)(0,20){5}{\line(-1,0){7}}
    \put(50,10){\circle*{3}}
    \put(50,50){\circle*{3}}
    \put(70,70){\circle*{3}}
    \put(90,30){\circle*{3}}
    \put(90,90){\circle*{3}}
    \put(150,10){\line(1,0){90}}
    \put(150,10){\line(0,1){90}}
    \multiput(150,10)(20,0){5}{\line(0,-1){7}}
    \multiput(150,10)(0,20){5}{\line(-1,0){7}}
    \put(150,10){\circle{3}}
    \put(150,30){\circle{3}}
    \put(170,10){\circle{3}}
    \put(170,30){\circle{3}}
     \put(190,10){\circle{3}}
    \multiput(150,50)(0,20){3}{\multiput(0,0)(20,0){2}{\circle*{3}}}
    \multiput(190,30)(0,20){4}{\circle*{3}}
    \multiput(210,10)(0,20){5}{\multiput(0,0)(20,0){2}{\circle*{3}}}
  \end{picture}
\caption{The elements of $A^\prime$ and the exponents of $I(A^\prime)$}
\label{ex2}
\end{figure}
\end{center}

Figure \ref{add2} displays the way in which $D(A^\prime)$ is built up from the three blocks 
$D(p_{1}^{-1}(1)\cap A^\prime)$,  $D(p_{1}^{-1}(2)\cap A^\prime)$ and $D(p_{1}^{-1}(3)\cap A^\prime)$. 
Note that here we do not simply stick the three blocks next to each other, as in the first example.

\begin{center}
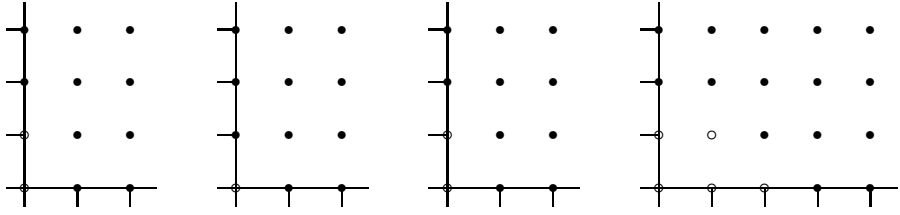
\begin{figure}[ht]
  \begin{picture}(390,90)
    \put(30,10){\line(1,0){50}}
    \put(30,10){\line(0,1){70}}
    \multiput(30,10)(20,0){3}{\line(0,-1){7}}
    \multiput(30,10)(0,20){4}{\line(-1,0){7}}
    \multiput(30,10)(0,20){2}{\circle{3}}
    \multiput(30,50)(0,20){2}{\circle*{3}}
    \multiput(50,10)(20,0){2}{\multiput(0,0)(0,20){4}{\circle*{3}}}
    \put(110,10){\line(1,0){50}}
    \put(110,10){\line(0,1){70}}
    \multiput(110,10)(20,0){3}{\line(0,-1){7}}
    \multiput(110,10)(0,20){4}{\line(-1,0){7}}
    \multiput(110,10)(20,0){3}{\line(0,-1){7}}
    \multiput(110,10)(0,20){4}{\line(-1,0){7}}
    \put(110,10){\circle{3}}
    \multiput(110,30)(0,20){3}{\circle*{3}}
    \multiput(130,10)(20,0){2}{\multiput(0,0)(0,20){4}{\circle*{3}}}
    \put(190,10){\line(1,0){50}}
    \put(190,10){\line(0,1){70}}
    \multiput(190,10)(20,0){3}{\line(0,-1){7}}
    \multiput(190,10)(0,20){4}{\line(-1,0){7}}
    \multiput(190,10)(20,0){3}{\line(0,-1){7}}
    \multiput(190,10)(0,20){4}{\line(-1,0){7}}
    \multiput(190,10)(0,20){2}{\circle{3}}
    \multiput(190,50)(0,20){2}{\circle*{3}}
    \multiput(210,10)(20,0){2}{\multiput(0,0)(0,20){4}{\circle*{3}}}
    \put(270,10){\line(1,0){90}}
    \put(270,10){\line(0,1){70}}
    \multiput(270,10)(20,0){5}{\line(0,-1){7}}
    \multiput(270,10)(0,20){4}{\line(-1,0){7}}
    \multiput(270,10)(20,0){2}{\multiput(0,0)(0,20){2}{\circle{3}}}
    \multiput(270,50)(20,0){2}{\multiput(0,0)(0,20){2}{\circle*{3}}}
    \put(310,10){\circle{3}}
    \multiput(310,30)(0,20){3}{\circle*{3}}
    \multiput(330,10)(20,0){2}{\multiput(0,0)(0,20){4}{\circle*{3}}}
  \end{picture}
\caption{$D(p_{1}^{-1}(1)\cap A^\prime)$,  $D(p_{1}^{-1}(2)\cap A^\prime)$ 
and $D(p_{1}^{-1}(3)\cap A^\prime)$ together form $D(A^\prime)$}
\label{add2}
\end{figure}
\end{center}

The ideas here presented can be generalised to arbitrary dimension $n$ and to arbitrary $A\subseteq\A^n$. 
We can sum the ideas up as follows.
\begin{itemize}
  \item We construct the Gr\"obner basis of $I(A)$ by induction over $n$. 
  \item The set $D(A)$ of those elements of $\N_{0}^n$ which do not occur as exponents of leading terms of elements of $I(A)$ 
    is built up from the sets $D(p_{1}^{-1}(a_{1})\cap A)\subseteq\N_{0}^{n-1}$ (analogous definition), 
    where $a_{1}$ runs through $p_{1}(A)$. (We will explain the way this is done in Sections \ref{addition} and \ref{gamma}.) 
  \item Assuming the induction hypothesis to hold true, we construct polynomials whose leading terms have exponents in
  $\N_{0}^{n-1}-D(p_{1}^{-1}(a_{1})\cap A)$, whose nonleading terms have exponents in 
  $D(p_{1}^{-1}(a_{1})\cap A)$, and which vanish on $p_{1}^{-1}(a_{1})\cap A$. 
  (This will be the content of the Corollary to Theorem \ref{thm} in Section \ref{it}.)
  \item These polynomials, along with one-dimensional interpolation, yield a collection of elements of $I(A)$. 
  (This collection will be constructed in Section \ref{almost}.) It does not form a Gr\"obner basis of $I(A)$, 
  unlike suitable linear combinations of elements of this collection. (This will be shown in the course of the proof of 
  Theorem \ref{thm} in Section \ref{it}.)
\end{itemize}

Finally, in Section \ref{comparison}, we will compare our construction of the Gr\"obner basis of $I(A)$ 
to the original method -- namely, the Buchberger--M\"oller algorithm.

\section{Notation}\label{notation}

We frequently use the projections
\begin{equation*}
  \begin{split}
    p_{i}:A&\to k\\
    (a_{1},\ldots,a_{n})&\mapsto a_{i}\text{ and}\\
    \widehat{p}^{i}:A&\to k^{n-1}\\
    (a_{1},\ldots,a_{n})&\mapsto (a_{1},\ldots,a_{i-1},a_{i+1},\ldots,a_{n})
  \end{split}
\end{equation*}
and will always write $\widehat{a}^{i}=\widehat{p}^{i}(a)$ for brevity's sake. 

As suggested in Section \ref{intro}, we frequently shift beween the use of monomials $X^\beta\in k[X]$ 
and the use of only their exponents $\beta\in\N_{0}^n$
via the equality $X_{i}=X^{e_{i}}$, where $e_{i}=(0,\ldots,0,1,0,\ldots,0)$, the $1$ situated at the $i$-th position.
In particular, we use the following collection of subsets of $\N_{0}^n$.
\begin{dfn}
  Let $\D_{n}$ be the set of all finite sets $D\subseteq\N_{0}^n$ such that whenever $d$ lies in $D$ and $d_{i}\neq0$, 
  then also $d-e_{i}$ lies in $D$. For $D\in\D_{n}$, we define its {\it limiting set} $E(D)$ to be the set of all 
  $\beta\in\N_{0}^n-D$ such that whenever $\beta_{i}\neq0$, then $\beta-e_{i}\in D$. 
\end{dfn}
Other characterisations of $E(D)$ are the following: $E(D)$ is the minimal subset $M\subseteq\N_{0}^n$ which generates
$\N_{0}^n-D$ as an additive submonoid of $\N_{0}^n$, or else: $E(D)$ is the minimal subset $M\subseteq\N_{0}^n$ 
satisfying $\cup_{\beta\in M}(\beta+\N_{0}^n)=\N_{0}^n-D$. 

\begin{center}
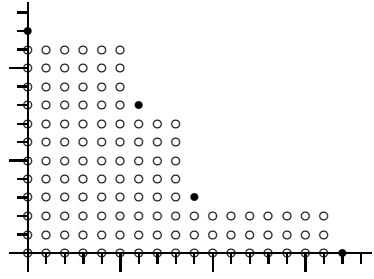
\begin{figure}[ht]
  \begin{picture}(160,115)
    \put(30,10){\line(1,0){130}}
    \put(30,10){\line(0,1){95}}
      \multiput(30,10)(7,0){19}{\line(0,-1){4}}
      \multiput(30,10)(0,7){14}{\line(-1,0){4}}
    \multiput(30,10)(35,0){4}{\line(0,-1){7}}
    \multiput(30,10)(0,35){3}{\line(-1,0){7}}  
    \multiput(30,10)(7,0){6}{\multiput(0,0)(0,7){12}{\circle{3}}}
    \multiput(72,10)(7,0){3}{\multiput(0,0)(0,7){8}{\circle{3}}}
    \multiput(93,10)(7,0){8}{\multiput(0,0)(0,7){3}{\circle{3}}}
    \put(30,94){\circle*{3}}
    \put(72,66){\circle*{3}}
    \put(93,31){\circle*{3}}
    \put(149,10){\circle*{3}}
  \end{picture}
\caption{An element of $\D_{2}$ and its limiting set}
\label{dee}
\end{figure}
\end{center}

In Figure \ref{dee}, the elements of some $D\in\D_{2}$ are marked with blank circles 
and the elements of $E(D)$ are marked with solid circles.

We can embed $\D_{n}$ into $\D_{n+1}$ by mapping each $d=(d_{1},\ldots,d_{n})\in D$ to $d=(d_{1},\ldots,d_{n},0)$
and conversely understand some $D\in\D_{n}$ such that $p_{n}(D)=\{0\}$ to lie in $\D_{n-1}$. This identification 
will become particularly important for our induction over $n$.

Also for $D\in\D_{n}$, the projections 
\begin{equation*}
  \begin{split}
    p_{i}:D&\to \N_{0}\text{ and}\\
    \widehat{p}^{i}:D&\to \N_{0}^{n-1}
  \end{split}
\end{equation*}
(defined by the same formulas as above) will be used frequently. As there is no danger of confusion, 
we do not use different names for the projections defined on $A$ and on some $D\in\D_{n}$. 

In fact, the only explicitly needed projections will be $p_{1}$ and $\widehat{p}^{1}$. Therefore, 
we write $p$ instead of $p_{1}$ and $\widehat{a}$ instead of $\widehat{p}^{1}(a)$. 
However, we will not replace $a_{1}=p(a)$ by any shorter notation.

A reader acquainted with the theory of Gr\"obner bases will of course immediately see the importance of $\D_{n}$ 
in the present context, but some words about the interpretation of $\D_{n}$ in terms of ideals of $k[X]$ may still be in order. 
Take a $D\subseteq\N_{0}^n$ and look at its complement $C=\N_{0}^n-D$. 
Then $D$ lies in $\D_{n}$ if, and only if, for all $c\in\N_{0}^n$ and for all $i\in\{1,\ldots,n\}$, 
$c\in C$ implies $c+e_{i}\in C$. 
Therefore, we can understand $C$ to be the set of exponents of leading terms of some ideal $J$ of $k[X]$. 
Let, for example, $J$ be the ideal generated by $X^\beta$, where $\beta$ runs through $E(D)$. 
Yet $J$ may also be assumed to be the ideal generated by a set of polynomials $f_{\beta}$, where $\beta$ runs through $E(D)$, 
such that the leading term of $f_{\beta}$ is $X^\beta$ and the exponents of all nonleading terms of $f_{\beta}$
lie in $D$. Conversely, a collection $(f_{\beta})_{\beta\in E(D)}$ of elements of an ideal $J$ 
of $k[X]$ is a Gr\"obner basis of $J$ precisely when for all $\beta$, the leading term of $f_{\beta}$ is $X^\beta$
and the exponents of all nonleading terms of $f_{\beta}$ lie in $D$. Hence for all ideals $J$, 
such a generating system exists. Furthermore, the finiteness of $D$ is equivalent to $k[X]/J$ 
being a finite dimensional $k$-vector space (or, $J$ being a zero-dimensional ideal). 
Therefore, $\D_{n}$ is precisely the set of those subsets of $\N_{0}^n$ 
that occur as exponents of non-leading terms of elements of some ideal $J$ of $k[X]$ such that $k[X]/J$ 
is a finite dimensional $k$-vector space.

Note that $k[X]/I(A)$ is a finite dimensional $k$-vector space by the Chinese Remainder Theorem. 
Therefore, the set $D(A)$ of exponents of nonleading terms of $I(A)$ lies in $\D_{n}$. 
Surprisingly, the apparition of $D(A)$ can be percieved simply from looking at the relative position of the elements of $A$.
The following two sections will deal with this.
 
\section{An addition map on $\D_{n}$}\label{addition}

The one cornerstone of our method is the following operation on $\D_{n}$.

\begin{dfn}
  For $D$ and $D^\prime$ in $\D_{n}$, define $D+D^\prime$ to be the set of all 
  $d\in\N_{0}^n$ such that $\widehat{d}\in\widehat{p}(D)\cup\widehat{p}(D^\prime)$ and 
  $d_{1}<\#\widehat{p}^{-1}(\widehat{d})\cap D+\#\widehat{p}^{-1}(\widehat{d})\cap D^\prime$.
\end{dfn}

To get a visual impression of what $+$ does, look at the example shown in Figure \ref{plus}. What is depicted there generalises
to arbitrary $D$ and $D^\prime$ in arbitrary dimension $n$ and can be described as follows. 
Draw a coordinate system of $\N_{0}^n$ and insert $D$. Place a translate of $D^\prime$ somewhere on the $1$-axis. 
The translate has to be sufficiently far out, so that $D$ and the translate of $D^\prime$ do not intersect. 
Then take the elements of the translate of $D^\prime$ and drop them down along the $1$-axis 
until they lie on top of an element of $D$, just as in the popular game {\it Connect4}. 
The result is $D+D^\prime$. 

\begin{center}
\begin{figure}[ht]
  \begin{picture}(336,115)
    \put(30,10){\line(1,0){130}}
    \put(30,10){\line(0,1){95}}
    \multiput(30,10)(7,0){19}{\line(0,-1){4}}
    \multiput(30,10)(0,7){14}{\line(-1,0){4}}
    \multiput(30,10)(35,0){4}{\line(0,-1){7}}
    \multiput(30,10)(0,35){3}{\line(-1,0){7}}  
    \multiput(30,10)(7,0){3}{\multiput(0,0)(0,7){13}{\circle{3}}}
    \multiput(51,10)(7,0){5}{\multiput(0,0)(0,7){3}{\circle{3}}}
    \put(51,84){$D$}
    \multiput(100,10)(7,0){5}{\multiput(0,0)(0,7){10}{\circle{3}}}
    \multiput(135,10)(7,0){3}{\multiput(0,0)(0,7){6}{\circle{3}}}
    \put(100,84){$(10,0)+D^\prime$}
    \put(190,10){\line(1,0){116}}
    \put(190,10){\line(0,1){95}}
    \multiput(190,10)(7,0){17}{\line(0,-1){4}}
    \multiput(190,10)(0,7){14}{\line(-1,0){4}}
    \multiput(190,10)(35,0){4}{\line(0,-1){7}}
    \multiput(190,10)(0,35){3}{\line(-1,0){7}} 
    \multiput(190,10)(7,0){3}{\multiput(0,0)(0,7){13}{\circle{3}}}
    \multiput(211,10)(7,0){5}{\multiput(0,0)(0,7){10}{\circle{3}}}
    \multiput(246,10)(7,0){3}{\multiput(0,0)(0,7){6}{\circle{3}}}
    \multiput(267,10)(7,0){5}{\multiput(0,0)(0,7){3}{\circle{3}}}
    \put(239,84){$D+D^\prime$}
   \end{picture}
\caption{Addition on $\D_{2}$}
\label{plus}
\end{figure}
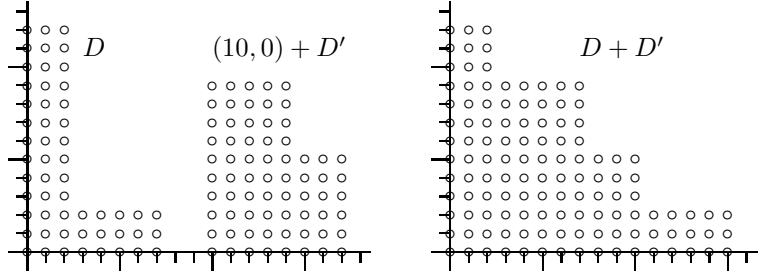
\end{center}

We will make use of the fact that 
$\widehat{p}(D+D^\prime)=\widehat{p}(D)\cup\widehat{p}(D^\prime)=\widehat{p}(D\cup D^\prime)$ 
and that $\#\widehat{p}^{-1}(\widehat{d})\cap(D+D^\prime)=\#\widehat{p}^{-1}(\widehat{d})\cap D
+\#\widehat{p}^{-1}(\widehat{d})\cap D^\prime$ for all $\widehat{d}\in\N_{0}^{n-1}$. 
Both being immediate consequences of the definition.  

\begin{lmm}\label{crit}
  Let $D\in\D_{n}$. Then for all $d\in D$, 
  $d_{1}<\#\widehat{p}^{-1}(\widehat{d})\cap D$. 
\end{lmm}

\begin{proof}
  If $d\in D$, then also all $d-\ell e_{1}\in D$ for all $\ell\in\{0,\ldots,d_{1}\}$.
\end{proof}

Of course, an analogous result holds true for all $i\in\{1,\ldots,n\}$, not only for $i=1$. 
Note that this gives the following characterisation of the limiting set of $D$: 
$\alpha\in E(D)$ if and only if
$\alpha_{i}=\#(\widehat{p}^{i})^{-1}(\widehat{d}^{i})\cap D$ for all $i\in\{1,\ldots,n\}$.

\begin{lmm}\label{lowerdim}
  If $D\in\D_{n}$, then $\widehat{p}(D)\in\D_{n-1}$.
\end{lmm}

\begin{proof}
  Take $d\in D$ with $d_{i}\neq0$. Then $d-e_{i}\in D$, hence
  $\widehat{d}-\widehat{e_{i}}=\widehat{d-e_{i}}\in\widehat{p}(D)$.
\end{proof}

\begin{lmm}\label{leq}
  Let $D\in\D_{n}$ and $d\in D$ such that $d_{i}\neq0$, where $i\neq 1$. Then 
  $\#\widehat{p}^{-1}(\widehat{d})\cap D\leq\#\widehat{p}^{-1}(\widehat{d-e_{i}})\cap D$.
\end{lmm}

\begin{proof}
  Define $\beta\in\N_{0}^n$ by setting $\widehat{\beta}=\widehat{d}$ and 
  $\beta_{1}=\#\widehat{p}^{-1}(\widehat{d})\cap D-1$. Then by Lemma \ref{crit}, 
  $\beta\in D$, hence also $\beta-e_{i}\in D$, and hence also $\beta-e_{i}-\ell e_{1}\in D$
  for all $\ell\in\{0,\ldots,\beta_{1}\}$. 
\end{proof}

\begin{pro}
  Let $D$, $D^\prime$ and $D^{\prime\prime}\in\D_{n}$. Then
  \begin{enumerate}
    \item[(a)] $D+D^\prime=D^\prime+D$,
    \item[(b)] $(D+D^\prime)+D^{\prime\prime}=D+(D^\prime+D^{\prime\prime})$,
    \item[(c)] $D+D^\prime\in\D_{n}$.
  \end{enumerate}
\end{pro}

\begin{proof}
  (a) This is clear.
  
  (b) The first set consists of those $d\in\N_{0}^n$ for which 
  $\widehat{d}\in\widehat{p}(D+D^\prime)\cup\widehat{p}(D^{\prime\prime})$ and 
  $d_{1}<\#\widehat{p}^{-1}(\widehat{d})\cap(D+D^\prime)+
  \#\widehat{p}^{-1}(\widehat{d})\cap D^{\prime\prime}$, which is the same as saying that
  $\widehat{d}\in\widehat{p}(D\cup D^\prime\cup D^{\prime\prime})$ and 
  $d_{1}<\#\widehat{p}^{-1}(\widehat{d})\cap D+
  \#\widehat{p}^{-1}(\widehat{d})\cap D^{\prime}+
  \#\widehat{p}^{-1}(\widehat{d})\cap D^{\prime\prime}$.
  In the same way, we can rewrite the conditions for $d$ to lie in the second set.
  
  (c) We have to show that for all $d\in D+D^\prime$ and for all $i\in\{1,\ldots,n\}$, 
  if $d_{i}\neq0$ then $d-e_{i}\in D+D^\prime$. 
  
  First let us look at $i=1$ and $d_{1}\neq0$. Then 
  $\widehat{d-e_{1}}=\widehat{d}\in\widehat{p}(D)\cup\widehat{p}(D^\prime)$
  and $(d-e_{1})_{1}=d_{1}-1<
  \#\widehat{p}^{-1}(\widehat{d})\cap D+\#\widehat{p}^{-1}(\widehat{d})\cap D^\prime$, 
  thus indeed $d-e_{1}\in D+D^\prime$.
  
  Now take $i\neq 1$ and $d_{i}\neq0$. 
  Due to symmetry, it suffices to consider the case where $\widehat{d}\in\widehat{p}(D)$. 
  By Lemma \ref{lowerdim}, we have $\widehat{d-e_{i}}=\widehat{d}-\widehat{e_{i}}\in
  \widehat{p}(D)\subseteq \widehat{p}(D)\cup\widehat{p}(D^\prime)$, and by Lemma \ref{leq}, we have
  $\#\widehat{p}^{-1}(\widehat{d})\cap D\leq\#\widehat{p}^{-1}(\widehat{d-e_{i}})\cap D$. 
  If we also have $\widehat{d}\in\widehat{p}(D^\prime)$, we analogously get 
  $\#\widehat{p}^{-1}(\widehat{d})\cap D^\prime
  \leq\#\widehat{p}^{-1}(\widehat{d-e_{i}})\cap D^\prime$. 
  And if $\widehat{d}\notin\widehat{p}(D^\prime)$, we simply have 
  $\#\widehat{p}^{-1}(\widehat{d})\cap D^\prime=0$, thus trivially also 
  $\#\widehat{p}^{-1}(\widehat{d})\cap D^\prime
  \leq\#\widehat{p}^{-1}(\widehat{d-e_{i}})\cap D^\prime$. 
  Subsumming inequalities, we get 
  $(d-e_{i})_{1}=d_{1}<\#\widehat{p}^{-1}(\widehat{d-e_{i}})\cap D+
  \#\widehat{p}^{-1}(\widehat{d-e_{i}})\cap D^\prime$, thus indeed $d-e_{i}\in D+D^\prime$.
\end{proof}

Therefore we can interpret $+$ as an addition map
\begin{equation*}
  \begin{split}
    +:\D_{n}\times\D_{n}&\to\D_{n}\\
    (D,D^\prime)&\mapsto D+D^\prime
  \end{split}
\end{equation*}
with the empty set as neutral element. In particular, given a finite family $(D_{b})_{b\in B}$ in $\D_{n}$, we can form the sum 
\begin{equation*}
    \sum_{b\in B}D_{b}\in\D_{n}\,.
\end{equation*}
This set can also be written as
\begin{equation*}
    \sum_{b\in B}D_{b}=\{d\in\N_{0}^n;
    \widehat{d}\in\underset{b\in B}{\cup}\widehat{p}(D_{b}),\,
    d_{1}<\sum_{b\in B}\#(\widehat{p})^{-1}(\widehat{d})\cap D_{b}\}\,.
\end{equation*}

\section{Assigning an element of $\D_{n}$ to the set of points}\label{gamma}

\begin{dfn}
  For $A\subseteq\A^n$ as above, we define $D(A)$ by induction over $n$ as follows.
  For $n=1$, we set $D(A)=\{0,\ldots,\#A-1\}$. 
  To pass from $n-1$ to $n$, we consider, for all $a_{1}\in p(A)$, the set $H(a_{1})=p^{-1}(a_{1})\cap A$. 
  We understand $H(a_{1})$ to be a subset of $\A^{n-1}$ via the projection map $\widehat{p}:H(a_{1})\to\A^{n-1}$.
  In this way, $D(H(a_{1}))$ is well-defined by the induction hypothesis. We set
  \begin{equation*}
    \begin{split}
      D(A)=\sum_{a_{1}\in p(A)}D(H(a_{1})).
    \end{split}
  \end{equation*}
\end{dfn}

Note that the induction might also be started at $n=0$ by defining $D(\text{Spec}\,k)=\{0\}$. 

Surely, the induction could also be disintangled by $D(H(a_{1}))$being written in terms of $D(H(a_{1},a_{2}))$, 
where $H(a_{1},a_{2})=(p_{1},p_{2})^{-1}(a_{1},a_{2})\cap A$, and so on. But this would not
make things more transparent.

Most is remarkable about this definition is that the set $D(A)$ we have just defined really is the set of exponents of 
nonleading terms of elements of $I(A)$. This will be elaborated upon in the Corollary to Theorem \ref{thm}. 

\section{A class of polynomials in $I(A)$}\label{almost}

With this we reach the second cornerstone of our method. For this, we take a closer look at the collection of $H(a_{1})$, 
where $a_{1}$ runs through $p(A)$, and on the respectice building blocks $D(H(a_{1}))$ of $D(A)$. 
As before, we understand $H(a_{1})$ to be a subset of 
$\A^{n-1}=\text{Spec}\,k[\widehat{X}]$, where $\widehat{X}=(X_{2},\ldots,X_{n})$. 

\begin{ass}
  We assume that for all $a_{1}\in p(A)$, the following holds true. 
  For all $\widehat{\beta}\in\N_{0}^{n-1}-D(H(a_{1}))$, there exists a polynomial 
  $f_{\widehat{\beta}}\in k[\widehat{X}]$ such that 
  \begin{enumerate}
    \item[(i)] the leading term of $f_{\widehat{\beta}}$ is $\widehat{X}^{\widehat{\beta}}$, 
    \item[(ii)] the exponents of all lower terms of $f_{\widehat{\beta}}$ lie in $\N_{0}^{n-1}-D(H(a_{1}))$, and
    \item[(iii)] $f_{\widehat{\beta}}(a)=0$ for all $a\in H(a_{1})$.
  \end{enumerate}
\end{ass}

$\beta\in E(D(A))$ given, let us split the set $p(A)$ into two components
\begin{equation*}
  \begin{split}
    S(\beta)&=\{a_{1}\in p(A);\widehat{\beta}\in D(H(a_{1}))\}\text { and}\\
    T(\beta)&=p(A)-S(\beta)\,.
  \end{split}
\end{equation*}
According to the above assumption, there is a polynomial $f_{\widehat{\beta},a_{1}}\in k[\widehat{X}]$ 
for all $a_{1}\in T(\beta)$ such that (i)--(iii) hold. Write this polynomial as
\begin{equation*}
  f_{\widehat{\beta},a_{1}}=\widehat{X}^{\widehat{\beta}}
  +\sum_{\widehat{\gamma}\in\widehat{G}_{a_{1}}}
  c_{\widehat{\beta},a_{1},\widehat{\gamma}}\widehat{X}^{\widehat{\gamma}}\,,
\end{equation*}
where $\widehat{\gamma}$ runs through the set 
\begin{equation*}
  \widehat{G}_{a_{1}}=\{\widehat{\gamma}\in D(H(a_{1}));\widehat{\gamma}<\widehat{\beta}\}\,.
\end{equation*}
We can even let $\widehat{\gamma}$ run through the bigger set 
\begin{equation*}
  \widehat{G}=\underset{a_{1}\in T(\beta)}{\cup}\widehat{G}_{a_{1}}
\end{equation*}
by simply setting $c_{\widehat{\beta},a_{1},\widehat{\gamma}}=0$ whenever 
$\widehat{\gamma}\in\widehat{G}-\widehat{G}_{a_{1}}$.
Next, we define
\begin{equation*}
  \theta_{\widehat{\beta}}=\widehat{X}^{\widehat{\beta}}+\sum_{a_{1}\in T(\beta)}
  \sum_{\widehat{\gamma}\in\widehat{G}}
  \chi(T(\beta),a_{1})c_{\widehat{\beta},a_{1},\widehat{\gamma}}\widehat{X}^{\widehat{\gamma}}\,,
\end{equation*}
where $\chi(T(\beta),a_{1})\in k[X_{1}]$ is the characteristic polynomial of $a_{1}\in T(\beta)$, i.e.,
\begin{equation*}
  \chi(T(\beta),a_{1})=\prod_{b_{1}\in T(\beta)-\{a_{1}\}}\frac{X_{1}-b_{1}}{a_{1}-b_{1}}\,.
\end{equation*}
Finally, we define
\begin{equation*}
  \phi_{\beta}=\prod_{a_{1}\in S(\beta)}(X_{1}-a_{1})\theta_{\beta}\,.
\end{equation*}
Let us state some properties of this polynomial.
\begin{itemize}
  \item The leading term of $\theta_{\widehat{\beta}}$ is $\widehat{X}^{\widehat{\beta}}$.
  Since $\beta\in E(D(A))$, Lemma \ref{crit} says that $\beta_{1}=\#S(\beta)$. 
  Therefore, the leading term of $\phi_{\beta}$ is $X^\beta$. 
  \item The exponents of all nonleading terms of $\theta_{\beta}$ lie in
  $\{0,\ldots,\#T(\beta)-1\}\times\widehat{G}$. Therefore, the exponents of the nonleading terms of
  $\phi_{\beta}$ lie in the union of $\{0,\ldots,\#p(A)-1\}\times\widehat{G}$ and 
  $\{0,\ldots,\#S(\beta)\}\times\{\widehat{\beta}\}$. 
  \item $\phi_{\beta}(a)=0$ for all $a\in A$. In fact, if $a_{1}\in S(\beta)$, this is most obvious. Else, 
  \begin{equation*}
    \theta_{\beta}(a)=\widehat{a}^{\widehat{\beta}}+\sum_{\widehat{\gamma}\in\widehat{G}}
    c_{\widehat{\beta},a_{1},\widehat{\gamma}}\widehat{a}^{\widehat{\gamma}}=
    \widehat{a}^{\widehat{\beta}}+\sum_{\widehat{\gamma}\in\widehat{G}_{a_{1}}}
    c_{\widehat{\beta},a_{1},\widehat{\gamma}}\widehat{a}^{\widehat{\gamma}}=
    f_{\widehat{\beta},a_{1}}(a)=0\,.
  \end{equation*}
\end{itemize}
So the polynomials $\phi_{\beta}$ really lie in $I(A)$ for all $\beta\in E(D(A))$.

\section{The main result}\label{it}

\begin{thm}\label{thm}
  Let $A\subseteq\A^n$ and $D(A)$ be as above, and let $\lambda\in E(D(A))$. 
  Then for all $\beta\in\cup_{\lambda^\prime\leq\lambda}(\lambda^\prime+\N_{0}^n)$, where 
  $\lambda^\prime$ runs through elements of $E(D(A))$, there is a polynomial $f_{\beta}\in k[X]$ such that 
  \begin{enumerate}
    \item[(i)] the leading term of $f_{\beta}$ is $X^\beta$, 
    \item[(ii)] the exponents of all lower terms of $f_{\beta}$ lie in 
      $\N_{0}^n-\cup_{\lambda^\prime\leq\lambda}(\lambda^\prime+\N_{0}^n)$, and
    \item[(iii)] $f_{\beta}(a)=0$ for all $a\in A$. 
  \end{enumerate}
\end{thm}

Before giving the proof, let us state and prove a corollary.

\begin{cor}
  For all $\beta\in \N_{0}^n-D(A)$, there is a unique $f_{\beta}\in k[X]$ such that 
  \begin{enumerate}
    \item[(i)] the leading term of $f_{\beta}$ is $X^\beta$, 
    \item[(ii)] the exponents of all lower terms of $f_{\beta}$ lie in $\N_{0}-D(A)$, and
    \item[(iii)] $f_{\beta}(a)=0$ for all $a\in A$. 
  \end{enumerate}
  In particular, the collection $f_{\beta}$, $\beta\in E(D(A))$ is a Gr\"obner basis of $I(A)$.
\end{cor}

\begin{proof}[Proof of Corollary]
  As for the existence of the polynomials $f_{\beta}$ as stated in the first part of the corollary, 
  let the Theorem be applied to the particular case where $\lambda$ is the maximal element of $E(D(A))$. 
  From this follows immediately that $(f_{\beta})_{\beta\in E(D(A))}$ is a Gr\"obner basis of $I(A)$. 
  In particular, the monomials $X^\gamma$, where $\gamma$ runs through $D(A)$, 
  are a basis of the $k$-vector space $k[X]/I(A)$. Now for the uniqueness of the polynomials $f_{\beta}$ 
  as stated in the first part of the corollary, assume that $g_{\beta}$ also satisfies properties (i)--(iii). 
  Then in particular $(f_{\beta}-g_{\beta})(a)=0$ for all $a\in A$, which means that $f_{\beta}-g_{\beta}$ lies in $I(A)$. 
  On the other hand, $f_{\beta}-g_{\beta}$ is an element of the $k$-span of $X^\gamma$, $\gamma\in D(A)$, 
  thus $f_{\beta}-g_{\beta}=0$. 
\end{proof} 

\begin{proof}[Proof of Theorem \ref{thm}]
  The proof will consist of 3 inductions, the outermost of which goes over $n\in\N$, the middle over $\lambda\in E(D(A))$ 
  and the innermost over $\beta\in\cup_{\lambda^\prime\leq\lambda}(\lambda^\prime+\N_{0}^n)$.
  
  So let us start with $n=1$. Here, the middle induction consists only of one the induction basis, since $E(D(A))=\{\#A\}$. 
  Therefore, we have to show that for all $\beta\geq\#A$, there is a unique polynomial $f_{\beta}\in k[X_{1}]$ with 
  properties (i)--(iii). 
  
  For $\beta=\#A$, take $f_{\beta}=\prod_{a\in A}(X-a)$. This polynomial clearly
  satisfies properties (i)--(iii). If the statement is shown for all $\beta^\prime<\beta$ in $\#A+\N_{0}$, 
  we define $f_{\beta}=Xf_{\beta-1}-c_{\#A}f_{\#A}$, where $c_{\#A}$ is the coefficient of 
  $X^{\#A-1}$ in $f_{\#A}$. This polynomial also satisfies properties (i)--(iii). 
  
  Thus, the statement is proved for $n=1$. The rest of the proof is the induction step from $n-1$ to $n$. So let $n>1$ be given.
  If the theorem is true for $n-1$, its corollary is true as well. Applying the corollary to the set 
  $H(a_{1})\subseteq\A^{n-1}$, we get precisely what we took for an assumption in the previous section. 
  Thus, we are given a collection $\phi_{\beta}$, for $\beta\in E(D(A))$, in $I(A)$, 
  as constructed in the previous section. We will presently make use of this collection.

  First, let $\lambda$ be the minimal element of $E(D(A))$. So we have to let $\beta$ run through all elements of
  $\lambda+\N_{0}^n$. 
  
  Let $\beta$ be the minimal element of this set, i.e., $\beta=\lambda$. By construction of $D(A)$, 
  we see that $\lambda=(\#p_{1}(A),0,\ldots,0)$. Thus we may take, analogously to what we have taken above, 
  $f_{\beta}=\prod_{a_{1}\in p_{1}(A)}(X_{1}-a_{1})$ and have properties (i)--(iii) satisfied. 
  
  Assume the statement is shown for all $\beta^\prime\in(\lambda+\N_{0}^n)$ 
  such that $\beta^\prime<\beta$. We show that the statement also holds true for $\beta$. 
  Since in this case $\beta^\prime$ is not equal to $\lambda$, there is an $i$ such that $\beta^\prime=\beta-e_{i}$ lies in
  $\lambda+\N_{0}^n$. Clearly $\beta^\prime<\beta$; therefore, the statement is true for $\beta^\prime$.   
  Consider the set $G$, which we define to be the set of all $\gamma\in(\lambda+\N_{0}^n)$
  such that $\gamma-e_{i}$ is the exponent of some nonleading term of $f_{\beta^\prime}$. 
  The statement is true also for all $\gamma\in G$, since if $\gamma-e_{i}$ 
  is the exponent of some nonleading term of $f_{\beta^\prime}$, then $\gamma-e_{i}<\beta^\prime=\beta-e_{i}$
  and therefore $\gamma<\beta$. Now we set
  \begin{equation}\label{G}
    f_{\beta}=X_{i}f_{\beta^\prime}-\sum_{\gamma\in G}c_{\gamma}f_{\gamma}\,,
  \end{equation}
  where $c_{\gamma}$ is the coefficient of $X^{\gamma-e_{i}}$ in $f_{\beta^\prime}$. 
  Again, properties (i)--(iii) are satisfied. 
  Thus the statement is shown for all $\beta\in(\lambda+\N_{0}^n)$. 
  
  Now we assume the statement is shown for all 
  $\beta\in\cup_{\lambda^\prime\leq\lambda^{\prime\prime}}(\lambda^{\prime}+\N_{0}^n)$, 
  where $\lambda^{\prime\prime}$ is the predecessor of $\lambda$ in $E(D(A))$. 
  We show that the statement is also true for all 
  $\beta\in\cup_{\lambda^\prime\leq\lambda}(\lambda^{\prime}+\N_{0}^n)$. 
  This will complete the proof of the theorem.
  
  First we note that the statement is true for all 
  $\beta\in\cup_{\lambda^\prime\leq\lambda}(\lambda^{\prime}+\N_{0}^n)$ such that $\beta<\lambda$. 
  In fact, $\beta$ even lies in 
  $\cup_{\lambda^\prime\leq\lambda^{\prime\prime}}(\lambda^{\prime\prime}+\N_{0}^n)$
  (otherwise $\beta\in(\lambda+\N_{0}^n)$ and therefore $\beta\geq\lambda$), which implies that there is an $f_{\beta}$ 
  with properties (i)--(iii). But in (ii), the exponents $\gamma$ of all lower terms of $f_{\beta}$ lie in 
  $\N_{0}-\cup_{\lambda^\prime\leq\lambda^{\prime\prime}}(\lambda^{\prime}+\N_{0}^n)$. 
  In fact, they even lie in $\N_{0}-\cup_{\lambda^\prime\leq\lambda}(\lambda^{\prime}+\N_{0}^n)$, 
  otherwise $\gamma\in(\lambda+\N_{0})$, thus $\beta<\lambda\leq\gamma$, which is a contradiction. 
  
  So we have to show the statement for all 
  $\beta\in\cup_{\lambda^\prime\leq\lambda}(\lambda^{\prime}+\N_{0}^n)$ such that $\beta\geq\lambda$. 
  The smallest such $\beta$ is $\beta=\lambda$. The polynomial $\phi_{\lambda}$ constructed in the previous section
  satisfies properties (i) and (iii) but not property (ii). To repair this, we have to get rid of all terms 
  of $\phi_{\lambda}$ whose exponents lie in 
  \begin{equation*}
    C=[0,\#p_{1}(A)-1]\times\{\widehat{\gamma}\in\cup_{T(\lambda)}D(H(a));
    \widehat{\gamma}<\widehat{\lambda}\}\,.
  \end{equation*}
  (Note that we do not have to get rid of those terms of $\phi_{\lambda}$ whose exponents lie in 
  $[0,\#S(\lambda)-1]\times\{\widehat{\lambda}\}$ since these lie in 
  $\N_{0}-\cup_{\lambda^\prime\leq\lambda}(\lambda^\prime+\N_{0})$, 
  as follows from the definition of $D(A)$.) Consider the set $G$, which we now define to be 
  $G=C\cap(\cup_{\lambda^\prime\leq\lambda}\lambda^\prime+\N_{0}^n)$. 
  The statement is shown for all $\gamma\in G$, since $\widehat{\gamma}<\widehat{\lambda}$ implies
  $\gamma<\lambda$ (in the lexicographic ordering). So the polynomial
  \begin{equation*}
    f_{\beta}=\phi_{\beta}-\sum_{\gamma\in G}c_{\gamma}f_{\gamma}\,,
  \end{equation*}
  where $c_{\gamma}$ is the coefficient of $X^{\gamma}$ in $\phi_{\lambda}$, is fine for properties (i)--(iii). 
  
  The last step is to assume that the statement is true for all $\beta^\prime<\beta$ in 
  $\cup_{\lambda^\prime\leq\lambda}(\lambda^{\prime}+\N_{0}^n)$ and to show that it is then also true for $\beta$. 
  Since we have already shown the statement for all $\beta$ equal to any of the $\lambda^{\prime}$ (which span 
  $\lambda^{\prime}+\N_{0}^n$), we now consider the complementary case. But in this case, there is an
  $i$ such that $\beta^\prime=\beta-e_{i}$ lies in
  $\cup_{\lambda^\prime\leq\lambda}(\lambda^{\prime}+\N_{0}^n)$. Therefore $\beta^\prime<\beta$, 
  thus the statement is true for $\beta^\prime$. The rest is analogous to what we did above.
  Define $G$ to be the set of all $\gamma\in\cup_{\lambda^\prime\leq\lambda}(\lambda^\prime+\N_{0}^n)$
  such that $\gamma-e_{i}$ is the exponent of some nonleading term of $f_{\beta^\prime}$. 
  Then for all $\gamma\in G$, the statement is true, since $\gamma-e_{i}<\beta^\prime=\beta-e_{i}$
  implies $\gamma<\beta$. The polynomial $f_{\beta}$, defined by the same formula as \eqref{G}, 
  satisfies properties (i)--(iii). And with this we are done.
\end{proof}

\section{Comparison with the Buchberger--M\"oller algorithm}\label{comparison}

Similarly to \cite{5}, let us give an informal description of how the Buchberger--M\"oller algorithm works. 

As already mentioned in Section \ref{intro}, the algorithm works not only for the lexicographic ordering on $k[X]$
but also for an arbitrary term ordering. In general, it is not clear what $D(A)$ looks like when $A$ is given. 
(Even in the case of lexicographic ordering, the shape of $D(A)$ has not been known before the present article.)
But since the exponents of the leading terms of the Gr\"obner basis of $I(A)$ are exactly the elements of $E(D(A))$,
one will have to determine $D(A)$ in one way or another. In the course of the Buchberger--M\"oller algorithm, 
this is done by considering one by one certain elements of $\N_{0}^n$ and deciding at each step whether or not 
the respective element belongs to $D(A)$. For so doing, one needs the following facts.

\begin{itemize}
  \item Assume we have found a subset $\Gamma$ of $D(A)$ which lies in $\D_{n}$. Take $\beta\in\N_{0}^n$ and define 
    $\Gamma^\prime=\Gamma\cup\{\beta\}$. Then $\beta$ does not lie in $D(A)$ if the rank of the matrix
    \begin{equation*}
      M(\Gamma^\prime)=(a^\gamma)_{\substack{\gamma\in\Gamma^\prime\\a\in A}}
    \end{equation*}
    is not maximal. 
  \item Conversely, that $M(\Gamma^\prime)$ is of maximal rank does not imply that $\beta\in D(A)$. 
    For this it is also necessary that $\beta$ be minimal amongst those elements of $\N_{0}^n-\Gamma$ that 
    might lie in $D(A)$. Call this set $B$. It consists of those elements of $E(\Gamma)$ 
    for which we have not yet checked the maximality of the rank of $M(\Gamma^\prime)$. 
\end{itemize}

Therefore, the algorithm for determining $D(A)$ goes as follows. 
\begin{itemize}
  \item Start with $\Gamma=\{0\}$ and $B=\{e_{1},\ldots,e_{n}\}$.
  \item When $\Gamma\in\D_{n}$ is given such that $\Gamma\subseteq D(A)$ and $B\subseteq E(D(A))$, 
    take $\beta$ to be the minimal element of $B$ and check whether the rank of $M(\Gamma^\prime)$ is maximal. 
    If this is the case, replace $\Gamma$ by $\Gamma\cup\{\beta\}$ and $B$ by $E(\Gamma\cup\{\beta\})$. 
    If not, just replace $B$ by $B-\{\beta\}$. 
  \item Proceed until $B=\emptyset$. In the end, $\Gamma=D(A)$. 
\end{itemize}

Once $D(A)$ is known, one computes a family of {\it separating polynomials} 
(the higher dimensional analogue of the characteristic polynomials in $k[X_{1}]$ we used). Consider the vector of polynomials 
\begin{equation*}
  (X^\gamma)=(X^\gamma)_{\gamma\in D(A)}\,.
\end{equation*}
Then the components of the vector of polynomials 
\begin{equation*}
  (\chi_{a})=(\chi_{a})_{a\in A}=M(\D(A))^{-1}(X^\gamma)
\end{equation*}
satisfy $\chi_{a}(a^\prime)=\delta_{a,a^\prime}$ for all $a$ and $a^\prime\in A$. 
Furthermore, the Gr\"obner basis of $I(A)$ is given by
\begin{equation*}
  f_{\beta}=X^\beta-\sum_{a\in A}X^\beta\chi_{a}\,,
\end{equation*}
where $\beta$ runs through all elements of $E(D(A))$. 

The method presented in this paper is fundamentally different in two ways. Firstly, we do not have to check whether any 
$\beta\in\N_{0}^n$ lies in $D(A)$, since we compute $D(A)$ by our inductive definition. 
Thus we save ourselves the trouble of computing the rank of $\#A$ matrices with $\leq\#A$ rows and $\#A$ columns. 
Secondly, we do not have to compute the inverse of the $\#A\times\#A$-matrix $M(D(A))$. 
However, during the course of Buchberger--M\"oller, one can compute the rank of the respective matrices in an 
iterative way, and in turn even successively compute the inverse of $M(D(A))$. 
This makes the algorithm far more effective, namely $\Oh(\#A^3)$. It is not hard to show that also our method is $\Oh(\#A^3)$. 
However, my personal judgement is that the virtue of our method lies in something else rather than in computational advantages: 
on the one hand, in the remarkable observation that we know what $D(A)$ will look like, 
and on the other hand in the insight we gain on what really makes up the elements of the Gr\"obner basis.

\section{Acknowledgements}

I wish to express gratitude to Leonhard Wieser, who introduced me to the Buchberger--M\"oller algorithm, 
and who gave me insight into his enlightning Master's thesis \cite{5}. I also thank Franz Pauer, advisor of Leonhard Wieser, 
who noticed a similarity of interests between Leo and me and therefore encouraged us to tell each other about our work. 
Special thanks go to Dominik Zeillinger (see \href{http://www.mathtics.doze.at/}{http://www.mathtics.doze.at/}),
who taught me how to transform the pictures from my head to a TeX file and to my brother Thomas Lederer, who helped 
me to transform my thoughts to comprehensible English.

\end{document}